\documentclass[10pt,letterpaper]{article}
\makeatletter
\usepackage{amscd}
\usepackage{amsmath}
\usepackage{latexsym}
\usepackage{amsfonts}
\usepackage{amssymb}
\usepackage{amsthm}
\usepackage{graphicx}

\newtheorem{thm}{\textbf Theorem}[section]

\newtheorem{rem}{\textbf Remark}[section]

\newcommand{\be}{\begin{eqnarray}}
\newcommand{\ee}{\end{eqnarray}}

\newcommand{\bes}{\begin{eqnarray*}}
\newcommand{\ees}{\end{eqnarray*}}

\begin{document}
\begin{titlepage}
\title{\bf Remarks on global regularity of 2D generalized MHD equations}
\author{ Baoquan Yuan\thanks{Corresponding Author: B. Yuan}\ and Linna Bai
       \\ School of Mathematics and Information Science,
       \\ Henan Polytechnic University,  Henan,  454000,  China.\\
        (bqyuan@hpu.edu.cn, blinna@163.com)
          }

\date{}
\end{titlepage}
\maketitle
\begin{abstract}
In this paper, we investigate the global regularity of 2D
generalized MHD equations, in which the dissipation term and
magnetic diffusion term are $\nu(-\Delta)^\alpha u$ and $\eta
(-\Delta)^\beta b$ respectively. Let $(u_{0}, b_{0})\in H^{s}$ with
$s\geq2$, it is showed that the smooth solution $(u(x,t),b(x,t))$ is
globally regular for the case $ 0\leq\alpha\leq\frac{1}{2},
\alpha+\beta > \frac{3}{2}.$

\end{abstract}

\vspace{.2in} {\bf Key words:}\quad Generalized MHD equations,
smooth solution, global regularity.



\section{Introduction}
\setcounter{equation}{0}
In this paper, we consider the following 2D generalized magnetohydrodynamic (GMHD) equations
 \be\label{1.1} \begin{cases}
{u_t} +\nu{\Lambda ^{2\alpha }}u+ u \cdot \nabla u =  - \nabla p + b \cdot \nabla b,\\
{b_t}+ \eta{\Lambda ^{2\beta }}b+ u \cdot \nabla b = b \cdot \nabla u,\\
\nabla\cdot u=\nabla\cdot b=0,\end{cases} \ee
  where $\alpha \ge 0$, $\beta\ge 0$, $\nu\geq0$ and $\eta\geq0$ are real parameters, and u is the velocity of the flow,
  b is the magnetic field, p is the scalar pressure, $\Lambda=(-\triangle)^{\frac{1}{2}}$ is defined in terms of Fourier transform by
  \bes \widehat{\Lambda f}(\xi)=|\xi|\widehat{ f(\xi)}. \ees
  If $\alpha=\beta=1$, (\ref{1.1}) is the viscous MHD equations, and the global wellposedness of classical solution is well-known \cite{S-T}.
  If $\nu=\eta=0$, (\ref{1.1}) is the invisid magnetohydrodynamic equations.

 We know that the 2D Euler equation is globally wellposed for
smooth initial data. But for the 2D invisid MHD equations, the
global wellposedness of classical solution is still a big open
problem. So the GMHD equations  has attracted much interest of many
mathematicians and has motivated a large number of research papers
concerning various generalizations and improvements \cite{T-Y-Z,
Wu1, Wu3, Wu2, Zhang, Zhou}. People pay attention to how the
parameters $\nu, \eta, \alpha, \beta$ influence the global
regularity of the GMHD equations. It is well-known that the
d-dimensional GMHD equations (\ref{1.1}) with $\nu>0$ and $\eta>0$
has a unique global classical solution for every initial data
$({u_0},{b_0}) \in {H^s}$ with $s\geq \max\{2\alpha, 2\beta\}$ if
$\alpha\geq\frac{1}{2}+\frac{d}{4}$ and
$\beta\geq\frac{1}{2}+\frac{d}{4}$ \cite{Wu1}. An improved result by
Wu \cite{Wu2} was established by reducing the requirement for
$\alpha$ and $\beta$ and the dissipation in (\ref{1.1}) by a
logarithmic factor. It is showed that the system is globally regular
as long as the following conditions
$\alpha\geq\frac{1}{2}+\frac{d}{4}$, $\beta>0$, $\alpha+\beta\geq
1+\frac{d}{2}$ are satisfied. As a special consequence, smooth
solutions of the 2D GMHD equations with $\alpha\geq1$, $\beta>0$,
$\alpha+\beta\geq2$ are global.

 However, for the 2D
incompressible MHD equations with partial dissipation, the global
regularity of the classical solutions is still a difficult problem.
In 2011, Cao and Wu \cite{C-W} showed an interesting result which
considered the 2D MHD equations of the form
 \be\label{1.2}
\begin{cases}
{u_t} + u \cdot \nabla u =  - p + {\nu _1}{u_{xx}} + {\nu _2}{u_{yy}} + b \cdot \nabla b,\\
{b_t} + u \cdot \nabla b = {\eta _1}{b_{xx}} + {\eta _2}{b_{yy}} + b \cdot \nabla u,\\
\nabla\cdot u=\nabla\cdot b=0,
\end{cases}
 \ee
they stated that the classical solutions of the equations
(\ref{1.2}) with either $\nu_{1}=0, \nu_{2}=\nu>0, \eta_{1}=\eta>0$
and $\eta_{2}=0$ or $\nu_{1}=\nu>0, \nu_{2}=0, \eta_{1}=0$ and
$\eta_{2}=\eta>0$ are globally existed for all time. If
$\nu_{1}=\nu_{2}=0$ and $\eta_{1}=\eta_{2}>0$, the MHD equations
(\ref{1.2}) has a global $H^{1}$ weak solution \cite{C-W, L-Z}. But
the existence of global classical solution is an open problem. When
$\eta_{1}=\eta_{2}=0$ and $\nu_{1}=\nu_{2}>0$, it is also unknown
for the existence of global classical solutions.

Recently, Tran, Yu and Zhai  \cite{T-Y-Z} obtained the global
regularity of 2D GMHD equations (\ref{1.2}) for the following three
cases: (1) $\alpha\geq\frac{1}{2}$, $\beta\geq 1$; (2)
$0\leq\alpha<\frac{1}{2}$, $2\alpha+\beta>2$; (3) $\alpha\geq2$,
$\beta=0$. Combining them with the result of \cite{Wu2}, we know
that if $\alpha+\beta\geq2$, (\ref{1.1}) with $\nu>0$ and $\eta>0$
possesses a global smooth solution. Note that in this case, the end
point $\alpha=0$ ($\nu=0$) and $\beta=2$ is not included and it
cannot ensure the global regularity for the system (\ref{1.1}).

Motivated by Tran, Yu and Zhai  \cite{T-Y-Z},  we carried on a
thorough investigation on whether the smooth solutions are global in
the case $\alpha=0$ and $\beta=2$ for 2D GMHD equations. In fact,
the system (\ref{1.1}) has a global classical solution for this
case. What is more, we find that when $\alpha=0$, the condition
$\beta=\alpha+\beta\geq2$ can be reduced to $\beta> \frac{3}{2}$.
When $0<\alpha\leq \frac{1}{2}$, We also conclude that  the system
is globally regular provided that $\alpha$ and $\beta$ satisfy the
relation $\alpha+\beta>\frac{3}{2}$.

To this end, we state our regularity criteria as follows.


 \begin{thm}\label{thm1}
  Consider the GMHD equations (\ref{1.1}) in  2D case. Assume $\left( {{u_0},{b_0}} \right) \in {H^s}$ with $s\geq2$.
  Then the system is globally regular for $\alpha$ and $\beta$ satisfying  $0\leq \alpha\leq\frac{1}{2}$,  $\alpha+\beta > \frac{3}{2}$.
\end{thm}
\begin{rem}
 In the special case $\alpha= \frac{1}{2}$, $\beta= 1$, reference \cite{T-Y-Z} showed that the equation (\ref{1.1}) is globally regular.
 However, the global regularity of (\ref{1.1}) with  $0\leq\alpha\leq \frac{1}{2}$, $\alpha+\beta=\frac{3}{2}$ is still a difficult problem.
 \end{rem}

\begin{rem} To simplify the presentation, we will set $\nu=\eta=1$.
It is a standard exercise to adjust various constants to accommodate
other values of $\nu$, $\eta$, as long as both are positive.
\end{rem}


\section{Proof of the main result}
\setcounter{equation}{0}

 In this section, we shall prove Theorem \ref{thm1}. The key idea here is to apply the standard $L^{2}-$ ennergy estimates to carry out the
 $H^{1}$, $H^{2}$ and higher estimates.


\subsection{ $L^{2}$ and $H^{1}$-energy estimates} We consider
the 2D GMHD equations (\ref{1.1}) with $\alpha\geq0$ and $\beta \geq
1$. it is easy to get the standard  $L^{2}$-energy estimate.
Multiplying the first two equations of (\ref{1.1}) by u and b,
respectively, integrating and adding the resulting equations
together it follows that \be \label{2.1} {\left\| u \right\|_2^2 +
\left\| b \right\|_2^2} + 2\int_0^t {\left\| {{\Lambda ^\alpha }u}
\right\|} _2^2\mbox{d}s + 2\int_0^t {\left\| {{\Lambda ^\beta }b}
\right\|} _2^2\mbox{d}s = \left\| {{u_0}} \right\|_2^2 + \left\|
{{b_0}} \right\|_2^2,\ee where we have used the incompressibility
condition $\nabla  \cdot u = \nabla  \cdot b = 0$.

As $\beta \geq 1$,  we can easily get
 \bes  b \in
{L^2}(0,T;{H^\beta(\mathbb{\mathbb{R}}^{2})})   \Rightarrow  \nabla
b \in {L^2}(0,T;{L^2(\mathbb{\mathbb{R}}^{2})}).
 \ees
 Let $\omega  =
\nabla  \times u =  - {\partial _2}{u_1} + {\partial _1}{u_2}$ be
the vorticity and $j = \nabla  \times b =  - {\partial _2}{b_1} +
{\partial _1}{b_2}$ be the current density.  Applying $\nabla\times$
to the first two equations of (\ref{1.1}) we obtain the governing
equations.
 \be
\begin{cases}{\omega _t} + u \cdot \nabla \omega  = b \cdot \nabla j  - {\Lambda ^{2\alpha }}\omega,\\
{j_t} + u \cdot \nabla j = b \cdot \nabla \omega  + T\left( {\nabla u,\nabla b} \right) - {\Lambda ^{2\beta }}j.\label{2.3}
\end{cases}
 \ee
Here
 \bes
T\left( {\nabla u,\nabla b} \right) = 2{\partial _1}{b_1}\left( {{\partial _1}{u_2}
+ {\partial _2}{u_1}} \right) + 2{\partial _2}{u_2}\left( {{\partial _1}{b_2} + {\partial _2}{b_1}} \right).
 \ees
Multiplying the two equations of (\ref{2.3}) by $\omega$ and j,  respectively,  integrating and applying the incompressibility condition we obtain
 \be
 \frac{1}{2}\frac{\mbox{d}}{{\mbox{d}t}}\int\limits_{{\mathbb{R}^2}} {({\omega ^2} + {j^2})} \mbox{d}x +
 \int\limits_{{\mathbb{R}^2}} {{{({\Lambda ^\alpha }\omega)}^2}} \mbox{d}x+ \int\limits_{{\mathbb{R}^2}} {{{({\Lambda ^\beta }j)}^2}} \mbox{d}x
 = \int\limits_{{\mathbb{R}^2}} {T\left( {\nabla u,\nabla b} \right)} j\mbox{d}x.\label{2.4}
 \ee
According to the Biot-Savart law, we have the representations
 \bes
 \frac{{\partial u}}{{\partial {x_k}}} = {R_k}(R \times \omega );   k = 1,2,
 \ees
  and
  \bes
  \frac{{\partial b}}{{\partial {x_k}}} = {R_k}(R \times j);    k = 1,2,
  \ees
where $R = ({R_1},{R_2})$, ${R_k} = - \frac{\partial {x_k}}{{( -
\Delta )^{ - \frac{1}{2}}}}$ denotes Riesz transformation. For
details about the Riesz transformation please refer to \cite{S-W}.
By the boundedness of Riesz operator R in $L^{p}$ space
$(1<p<\infty)$,  we arrive at
 \bes
{\left\| {\nabla u} \right\|_{{L^2}}} \le C{\left\| \omega
\right\|_{{L^2}}} \quad and \quad  {\left\| {\nabla b}
\right\|_{{L^4}}} \le C{\left\| j \right\|_{{L^4}}}.
 \ees
  Using H\"{o}lder and Young's inequalities one has
   \bes
\int\limits_{{\mathbb{R}^2}} {T\left( {\nabla u,\nabla b} \right)}
j\mbox{d}x &&\le C{\left\| {\nabla u} \right\|_{{L^2}}}{\left\|
{\nabla b} \right\|_{{L^4}}}{\left\| j \right\|_{{L^4}}} \le
C{\left\| \omega  \right\|_{{L^2}}}\left\| j \right\|_{_{{L^2}}}^{2
- \frac{1}{\beta }}\left\| {{\Lambda ^\beta }j}
\right\|_{_{{L^2}}}^{\frac{1}{\beta }}\\&&\le C(\varepsilon )\left\|
\omega  \right\|_{{L^2}}^{\frac{{2\beta }}{{2\beta  - 1}}}\left\| j
\right\|_{{L^2}}^2 + \varepsilon \left\| {{\Lambda ^\beta }j}
\right\|_{{L^2}}^2\le C(\varepsilon )(\left\| \omega
\right\|_{{L^2}}^2 + 1)\left\| j \right\|_{{L^2}}^2 + \varepsilon
\left\| {{\Lambda ^\beta }j} \right\|_{{L^2}}^2,
 \ees
 where we have used the following Gagliardo-Nirenberg inequality
  \bes {\left\| j
\right\|_{{L^4}}} \le C\left\| j \right\|_{{L^2}}^{1 -
\frac{1}{{2\beta }}}\left\| {{\Lambda ^\beta }j}
\right\|_{{L^2}}^{\frac{1}{{2\beta }}}.
 \ees
 Inserting the above estimate into (\ref{2.4}),  and taking $\varepsilon$ small enough so
that $\varepsilon<1$ we have
 \bes
 \frac{\mbox{d}}{\mbox{d}t}
\left(\left\| \omega  \right\|_{{L^2}}^2 + \left\| j
\right\|_{{L^2}}^2 \right)+\left\| {{\Lambda ^\alpha }\omega }
\right\|_{{L^2}}^2 + \left\| {{\Lambda ^\beta }j} \right\|_{{L^2}}^2
\leq C(\varepsilon)\left( \left\| \omega  \right\|_{{L^2}}^2 + 1
\right)\left\| j \right\|_{{L^2}}^2 .
 \ees
 Gronwall's inequality \cite[Appenddix B.j]{Evans} and $L^{2}$ energy estimate imply that
  \bes
  \left\| \omega  \right\|_{{L^2}}^2 + \left\| j \right\|_{{L^2}}^2
  +\int_0^t {\left\| {{\Lambda ^\alpha }\omega } \right\|_{{L^2}}^2\mbox{d}s
  + \int_0^t {\left\| {{\Lambda ^\beta }j} \right\|_{{L^2}}^2\mbox{d}s} }
  \le \left( {\left\| {{\omega _0}} \right\|_{{L^2}}^2
  + \left\| {{j_0}} \right\|_{{L^2}}^2} \right)
  \exp\left[ {\int_0^t {\left\| j \right\|_{{L^2}}^2} \mbox{d}s } \right]<\infty.
  \ees


 \subsection{ Higher estimates for $\alpha=0$}
 In this case we have $\beta > \frac{3}{2}$,  and the GMHD equations now read
 \be
 \label{u}
 \begin{cases}
{u_t} + u \cdot \nabla u =  - \nabla p + b \cdot \nabla b, \\
{b_t} + u \cdot \nabla b = b \cdot \nabla u - {\Lambda ^{2\beta }}b,\\
\nabla\cdot u=\nabla\cdot b=0.
 \end{cases}
 \ee
First of all, we estimate $b_{t}$.
Taking the inner product of the second equation of (\ref{u}) with $b_{t}$ and using H\"{o}lder and  Young's inequalities we obtain
 \bes
\left\| {{b_t}} \right\|_{{L^2}}^2 + \frac{1}{2}\frac{\mbox{d}}{{\mbox{d}t}}\left\| {{\Lambda ^\beta }b} \right\|_{{L^2}}^2
&&\le \int\limits_{{\mathbb{R}^2}} {\left| {u \cdot \nabla b \cdot {b_t}} \right|\mbox{d}x}
+ \int\limits_{{\mathbb{R}^2}} {\left| {b \cdot \nabla u \cdot {b_t}} \right|\mbox{d}x}\\
&&\le \frac{1}{2}\left\| {{b_t}} \right\|_{{L^2}}^2 + \frac{1}{2}\left( {\left\| u \right\|_{{L^4}}^2
\left\| {\nabla b} \right\|_{{L^4}}^2 + \left\| {\nabla u} \right\|_{{L^2}}^2\left\| b \right\|_{{L^\infty }}^2} \right).
 \ees
Application of the following Gagliardo-Nirenberg inequalities
 \bes
&&{\left\| f \right\|_{{L^4}}} \le C\left\| f \right\|_{{L^2}}^{\frac{1}{2}}\left\| {\nabla f} \right\|_{{L^2}}^{\frac{1}{2}},\\
&&{\left\| f \right\|_{{L^\infty }}} \le C\left\| f \right\|_{{L^2}}^{1 - \frac{1}{\beta }}\left\| {{\Lambda ^\beta }f} \right\|_{{L^2}}^{\frac{1}{\beta }},
 \ees
 yields that
\bes
 \left\| {{b_t}} \right\|_{{L^2}}^2 +\frac{\mbox{d}}{{\mbox{d}t}}\left\| {{\Lambda ^\beta }b} \right\|_{{L^2}}^2
 &&\le  C{\left\| u \right\|_{{L^2}}}{\left\| {\nabla u} \right\|_{{L^2}}}{\left\| {\nabla b} \right\|_{{L^2}}}{\left\| {\nabla j} \right\|_{{L^2}}}
 + C\left\| {\nabla u} \right\|_{{L^2}}^2\left\| b \right\|_{{L^2}}^{1 - \frac{1}{\beta }}\left\| {{\Lambda ^\beta }b} \right\|_{{L^2}}^{\frac{1}{\beta }}\\
 &&\le C{\left\| {\nabla j} \right\|_{{L^2}}} + C\left\| {{\Lambda ^\beta }b} \right\|_{{L^2}}^{\frac{1}{\beta }}.
  \ees
By the results of the $L^{2}$-energy estimate and $H^{1}$ estimate,  we deduce that
\be\label{2.2}
\left\| {{\Lambda ^\beta }b} \right\|_{{L^2}}^2 + \int_0^t {\left\| {{b_t}} \right\|_{{L^2}}^2\mbox{d}s}
\le \left\| {{\Lambda ^\beta }{b_0}} \right\|_{{L^2}}^2+C\int_0^t {\left\| {{\Lambda ^\beta }b} \right\|_{{L^2}}^{\frac{1}{\beta }}\mbox{d}s + }
C\int_0^t {{{\left\| {\nabla j} \right\|}_{{L^2}}}\mbox{d}s} < \infty .
 \ee
Now we go back to the equation ${b_t} + u \cdot \nabla b = b \cdot \nabla u - {\Lambda ^{2\beta }}b$,
and using the similar way with the estimate of $b_{t}$ we get
 \bes
\left\| {{\Lambda ^{2\beta }}b} \right\|_{{L^2}}^2\le \left\| {{b_t}} \right\|_{{L^2}}^2
+ \left\| {u \cdot \nabla b} \right\|_{{L^2}}^2
+ \left\| {b \cdot \nabla u} \right\|_{{L^2}}^2\le \left\| {{b_t}} \right\|_{{L^2}}^2
+ C{\left\| {\nabla j} \right\|_{{L^2}}} + C\left\| {{\Lambda ^\beta }b} \right\|_{{L^2}}^{\frac{1}{\beta }}.
 \ees
Recall that $j=\nabla\times b$,  one can deduce,  thanks to (\ref{2.2}),  that
 \be \label{bbb}
 &&\nonumber\int_0^t {{{\left\| {\nabla j} \right\|}_{{{\dot{H}}^{2\beta-2}}}^2} \mbox{d}s}
 \le \int_0^t {\left\| {{\Lambda ^{2\beta }}b} \right\|_{{L^2}}^2} \mbox{d}s \\
 && \le \int_0^t {\left\| {{b_t}} \right\|_{{L^2}}^2\mbox{d}s}  + C\int_0^t {{{\left\| {\nabla j} \right\|}_{{L^2}}}} \mbox{d}s
 + C\int_0^t {\left\| {{\Lambda ^\beta }b} \right\|_{{L^2}}^{\frac{1}{\beta }}\mbox{d}s}< \infty.
 \ee
 Since $\beta > \frac{3}{2}$,  by Sobolev embedding theorem,  it is easily to see
 \bes
 \nabla j \in {L^2}(0,T;{H^{2\beta  - 2}(\mathbb{R}^{2})})\hookrightarrow  {L^2}(0,T;{L^\infty (\mathbb{\mathbb{R}}^{2})}).
 \ees
Secondly,  we estimate $\omega$.  From the first equation of (\ref{u}),  we have the vorticity equation
$\omega_{t}+u\cdot\nabla\omega=b\cdot\nabla j$.
Multiplying both sides of it by $p{\left| \omega  \right|^{p - 2}}\omega $ and integrating both sides over $\mathbb{R}^{2}$,
it follows, by H\"{o}lder inequality,  that
 \bes
 \frac{\mbox{d}}{{\mbox{d}t}}\left\| \omega  \right\|_{{L^p}}^p
 + p{\int\limits_{{\mathbb{R}^2}} {u \cdot \nabla \omega \cdot \left| \omega  \right|} ^{p - 2}}\omega \mbox{d}x
 \le p{\left\| {b \cdot \nabla j} \right\|_{{L^p}}}\left\| \omega  \right\|_{{L^p}}^{p - 1}
 \leq p{\left\| b \right\|_{{L^\infty }}}{\left\| {\nabla j} \right\|_{{L^p}}}\left\| \omega  \right\|_{{L^p}}^{p - 1}.
 \ees
Noting that  $p{\int_{\mathbb{R}^2} {u \cdot \nabla \omega \cdot
\left| \omega  \right|} ^{p - 2}}\omega \mbox{d}x=0$. Now let $p \to
\infty $, we infer that
 \bes
 {\left\| \omega  \right\|_{{L^\infty }}} \le {\left\| {{\omega _0}} \right\|_{{L^\infty }}}
 + \int_0^t {{{\left\| b \right\|}_{{L^\infty }}}{{\left\| {\nabla j} \right\|}_{{L^\infty }}}\mbox{d}s}  < \infty.
  \ees
This leads to
 \bes
 \omega  \in {L^\infty }(0,T;{L^\infty(\mathbb{R}^{2}) }).
  \ees
 Lastly,  according to the classical BKM-type blow up criterion \cite{C-K-S} which is the MHD system stays regular beyond T
 provided that $\int_0^T ({{{\left\| \omega  \right\|}_{{L^\infty }}} + {{\left\| j \right\|}_{{L^\infty }}}})\mbox{d}t  < \infty $,
  the proof of the case $\alpha=0$ is thus completed.


\subsection{Higher estimates for $0<\alpha \leq \frac{1}{2}$,
$\alpha+\beta> \frac{3}{2}$}

In this case,  we can easily get $\beta>1$. Firstly,  we estimate
${{{\left\| \omega \right\|}_{{L^p}}}}$. Multiplying both sides of
the first equation of (\ref{2.3}) by $p{\left| \omega  \right|^{p -
2}}\omega $ and integrating both sides over $\mathbb{R}^{2}$,  it
follows that
 \bes
 \frac{\mbox{d}}{{\mbox{d}t}}\left\| \omega
\right\|_{{L^p}}^p + p\int\limits_{{\mathbb{R}^2}} {{\Lambda
^{2\alpha }}\omega  \cdot {{\left| \omega  \right|}^{p - 2}}\omega
\mbox{d}x}  \le p{\left\| {b \cdot \nabla j}
\right\|_{{L^p}}}\left\| \omega \right\|_{{L^p}}^{p - 1}.
 \ees
For the dissipation term,  we know by the property of Riesz
potential that $\int\limits_{{\mathbb{R}^2}} {{\Lambda ^{2\alpha
}}\omega \cdot {{\left| \omega  \right|}^{p - 2}}\omega
\mbox{d}x}\geq0 $. For the details on it see \cite{A-D}.
 Thus,  we have
  \be {\left\| \omega
\right\|_{{L^p}}} \le {\left\| {{\omega _0}} \right\|_{{L^p}}} +
\int_0^t {{{\left\| {b \cdot \nabla j} \right\|}_{{L^p}}}\mbox{d}s}.
\label{2.6}
 \ee
  By the Gagliardo-Nirenberg inequality,  one has the
following estimate
 \be \label{2.8}
 {\left\| {b \cdot \nabla j}
\right\|_{{L^p}}} \le C{\left\| b \right\|_{{L^\infty }}}{\left\|
{\nabla j} \right\|_{{L^p}}} \le C\left\| b
\right\|_{{L^2}}^{\frac{\beta }{{1 + \beta }}}\left\| {{\Lambda
^\beta }j} \right\|_{_{{L^2}}}^{\frac{1}{{1 + \beta }}}\left\| j
\right\|_{{L^2}}^{\frac{{2\beta  - 3}}{{2\beta  - 1}} +
\frac{2}{{(2\beta  - 1)p}}}\left\| {{\Lambda ^{2\beta  - 1}}j}
\right\|_{_{{L^2}}}^{\frac{{2(p - 1)}}{{(2\beta  - 1)p}}},
 \ee
where $p$ satisfies $p > \frac{1}{\alpha }$. So,  inserting
(\ref{2.8}) into (\ref{2.6}) and applying with $L^{2}$ and $H^{1}$
estimates and (\ref{bbb}),  it can be derived that
 \bes
 {\left\|
\omega  \right\|_{{L^p}}}&& \le {\left\| {{\omega _0}}
\right\|_{{L^p}}} + C\int_0^t {\left\| b
\right\|_{{L^2}}^{\frac{\beta }{{1 + \beta }}}\left\| {{\Lambda
^\beta }j} \right\|_{_{{L^2}}}^{\frac{1}{{1 + \beta }}}\left\| j
\right\|_{{L^2}}^{\frac{{2\beta  - 3}}{{2\beta  - 1}} +
\frac{2}{{(2\beta  - 1)p}}}\left\| {{\Lambda ^{2\beta  - 1}}j}
\right\|_{_{{L^2}}}^{\frac{{2(p - 1)}}{{(2\beta  - 1)p}}}\mbox{d}s}\\
&&\le {\left\| {{\omega _0}} \right\|_{{L^p}}} + C\int_0^t {\left\|
{{\Lambda ^\beta }j} \right\|_{_{{L^2}}}^{\frac{1}{{1 + \beta
}}}\left\| {{\Lambda ^{2\beta  - 1}}j}
\right\|_{_{{L^2}}}^{\frac{{2(p - 1)}}{{(2\beta  - 1)p}}}\mbox{d}s}\\
&&\le {\left\| {{\omega _0}} \right\|_{{L^p}}} + C\int_0^t {(\left\|
{{\Lambda ^\beta }j} \right\|_{_{{L^2}}}^2 + \left\| {{\Lambda
^{2\beta }}b} \right\|_{_{{L^2}}}^{\frac{{4(p - 1)(1 + \beta
)}}{{(2\beta  - 1)(2\beta  + 1)p}}})\mbox{d}s}  < \infty.
 \ees
Note that as long as $p>\frac{1}{\alpha}$,  we have $\frac{{4(p -
1)(1 + \beta )}}{{(2\beta  - 1)(2\beta  + 1)p}} \le 2$.

Secondly, we derive the estimates of  ${{{\left\| \omega
\right\|}_{{H^1}}}}$ and ${{{\left\| j  \right\|}_{{H^1}}}}$.  We
differentiate the equations (\ref{2.3}) with respect to $x_{i}$ over
$\mathbb{R}^{2}$,  then multiply the resulting equations by
${\partial _{{x_i}}}\omega $ and ${\partial _{{x_i}}}j$ for $i=1,
2$,  integrate with respect to $x$ over $\mathbb{R}^{2}$ and sum
them up.  It follows that
 \be \nonumber
&&\frac{1}{2}\frac{\mbox{d}}{{\mbox{d}t}}(\left\| {\nabla \omega }
\right\|_{{L^2}}^2 + \left\| {\nabla j} \right\|_{{L^2}}^2) +
\left\| {{\Lambda ^\alpha }\nabla \omega } \right\|_{{L^2}}^2 +
\left\| {{\Lambda ^\beta }\nabla j} \right\|_{{L^2}}^2\\ \nonumber
&&\leq \int {\left| {\nabla u} \right|} {\left| {\nabla \omega }
\right|^2}\mbox{d}x + \int {\left| {\nabla b} \right|} \left|
{\nabla j} \right|\left| {\nabla \omega } \right|\mbox{d}x  + \int
{\left| {\nabla u} \right|} {\left| {\nabla j} \right|^2}\mbox{d}x
\\ \nonumber &&+ \int {\left| {\nabla b} \right|} \left| {\nabla
\omega } \right|\left| {\nabla j} \right|\mbox{d}x + \int {\left|
{{\nabla ^2}u} \right|} \left| {\nabla b} \right|\left| {\nabla j}
\right|\mbox{d}x + \int {\left| {\nabla u} \right|} \left| {{\nabla
^2}b} \right|\left| {\nabla j} \right|\mbox{d}x \\   &&= {I_1} +
{I_2} + {I_3} + {I_4} + {I_5} + {I_6}.\label{2.7}
 \ee
It is easy to see that the estimates of ${I_4}$ and ${I_5}$ are the
same as ${I_2}$ while ${I_6}$ is the same as  ${I_3}$.  Therefore,
it suffices to estimate $I_{1}$, $I_{2}$, $I_{3}$.\\
H\"{o}lder, Young and Galiardo-Nirenberg inequalities together give
 \bes
 {I_1} \le {\left\| {\nabla u} \right\|_{{L^p}}}\left\| {\nabla
\omega } \right\|_{{L^{2q}}}^2 \le C{\left\| \omega
\right\|_{{L^p}}}\left\| {\nabla \omega }
\right\|_{{L^2}}^{\frac{{2(\alpha p - 1)}}{{\alpha p}}}\left\|
{{\Lambda ^\alpha }\nabla \omega }
\right\|_{{L^2}}^{\frac{2}{{\alpha p}}} \le C(\varepsilon)\left\|
{\nabla \omega } \right\|_{{L^2}}^2 + \varepsilon \left\| {{\Lambda
^\alpha }\nabla \omega } \right\|_{{L^2}}^2,
 \ees
where $p$ and $q$ satisfy $\frac{1}{p}+\frac{1}{q}=1$ and $p>\frac
 1\alpha$.\\
Arguing similarly as the estimate of $I_{1}$,  thanks to the $L^{2}$
and $H^{1}$  estimates, one has
 \bes
{I_2}&&\le {\left\| {\nabla b} \right\|_{{L^\infty}}}{\left\|
{\nabla j} \right\|_{{L^2}}}{\left\| {\nabla \omega }
\right\|_{{L^2}}} \le C\left\| \nabla b
\right\|_{{L^2}}^{1-\frac{1}{\beta}}\left\| {\Lambda^{\beta}\nabla
b} \right\|_{{L^2}}^{\frac{1}{\beta}}( \left\| {\nabla \omega}
\right\|_{{L^2}}^{2}+ \left\|\nabla j \right\|_{{L^2}}^{2})
\\&&\le C\|{\Lambda^{\beta}\nabla
b}\|_{{L^2}}^{\frac1\beta}(\|\nabla\omega\|_{L^{2}}^{2}+\|\nabla
j\|_{L^{2}}^{2})
 \ees
where  use has been made of the following Gagliardo-Nirenberg
inequality
 \bes
{\left\| {\nabla b} \right\|_{{L^\infty}}} \le C\left\| \nabla b
\right\|_{{L^2}}^{1-\frac{1}{\beta}}\left\| {{\Lambda ^\beta }\nabla
b} \right\|_{{L^2}}^{\frac{1}{{\beta }}}.
 \ees
The estimate of $I_{3}$ can also be obtained by H\"{o}lder, Young
and Sobolev embedding inequalities
 \bes
 {I_3}\le {\left\| {\nabla u} \right\|_{{L^2}}}\left\| {\nabla j}
\right\|_{{L^4}}^2 \le C{\left\| {\nabla u} \right\|_{{L^2}}}\left\|
{\nabla j} \right\|_{{L^2}}^{\frac{{2\beta  - 1}}{\beta }}\left\|
{{\Lambda ^\beta }\nabla j} \right\|_{{L^2}}^{\frac{1}{\beta }} \le
C(\varepsilon )\left\| {\nabla j} \right\|_{{L^2}}^2 + \varepsilon
\left\| {{\Lambda ^\beta }\nabla j} \right\|_{{L^2}}^2.
 \ees
Combining the above estimates into (\ref{2.7}),  and taking
$\varepsilon$ small enough we get
 \bes
\frac{1}{2}\frac{{\rm{d}}}{{{\rm{d}}t}}(\left\| {\nabla \omega }
\right\|_{{L^2}}^2 + \left\| {\nabla j} \right\|_{{L^2}}^2)+\left\|
{{\Lambda ^\alpha }\nabla \omega } \right\|_{{L^2}}^2 + \left\|
{{\Lambda ^\beta }\nabla j} \right\|_{{L^2}}^2 \le
C(\varepsilon)\|{\Lambda^{\beta}\nabla
b}\|_{{L^2}}^{\frac1\beta}(\left\| {\nabla \omega }
\right\|_{{L^2}}^2 + \left\| {\nabla j} \right\|_{{L^2}}^2).
 \ees
Gronwall's inequality and $H^{1}$ estimate imply that
 \bes
\left\| {\nabla \omega } \right\|_{{L^2}}^2 + \left\| {\nabla j}
\right\|_{{L^2}}^2 +\int_0^t {\left\| {{\Lambda ^\alpha }\nabla
\omega } \right\|_{{L^2}}^2\mbox{d}s + \int_0^t {\left\| {{\Lambda
^\beta }\nabla j} \right\|_{{L^2}}^2} } \mbox{d}s \le C(\left\|
{\nabla {\omega _0}} \right\|_{{L^2}}^2 + \left\| {\nabla {j_0}}
\right\|_{{L^2}}^2 ).
 \ees
Thus,we arrive at
 \bes
&&\omega  \in {L^\infty }(0,T;{H^1}({\mathbb{R}^2}))\cap {L^2 }(0,T;{H^{\alpha+1}}({\mathbb{R}^2})),\\
&& j  \in {L^\infty }(0,T;{H^1}({\mathbb{R}^2}))\cap {L^2
}(0,T;{H^{\beta+1}}({\mathbb{R}^2})).
 \ees
 In the end, by the embedding  relation $H^{s}(\mathbb{R}^{2})\hookrightarrow L^{\infty}(\mathbb{R}^{2})$ for $s>1$,
 we can get $\omega  \in {L^2 }(0,T;L^{\infty}({\mathbb{R}^2})) ,  j  \in {L^2 }(0,T;L^{\infty}({\mathbb{R}^2}))$,
 and combining the BKM-type blow-up criterion \cite{C-K-S},  this completes the proof.
 Obviously, the fact $H^{1}(\mathbb{R}^{2})\hookrightarrow \mbox{BMO}(\mathbb{R}^{2})$ and  the blow-up criterion \cite{Z-L} can also give the proof.  $\Box$ \\
 \vspace{0.4cm}

 \textbf{Acknowledgements} The research of B Yuan
was partially supported by the National Natural Science Foundation
of China (No. 11071057), Innovation Scientists and Technicians Troop
Construction Projects of Henan Province (No. 104100510015).


\end{document}